\documentclass{amsart}
\usepackage{amsmath}
\usepackage{amsfonts,amssymb}
\usepackage[latin1]{inputenc}

\newtheorem{theorem}{Theorem}
\newtheorem{corollary}[theorem]{Corollary}
\newtheorem{lemma}[theorem]{Lemma}
\newtheorem{definition}[theorem]{Definition}
\newtheorem{proposition}[theorem]{Proposition}

\newtheorem{remark}[theorem]{Remark}

\newcommand{\fx}{f(x_1,\ldots,x_n)}
\newcommand{\F}{F\langle X|G\rangle}
\newcommand{\FF}{F\langle X\rangle}
\newcommand{\Z}{\mathbb{Z}}

\newcommand{\zn}{\mathbb{Z}_n}

\newcommand{\U}{\mathcal{U}}
\newcommand{\support}{\text{\rm Supp}}

\begin{document}
\title[Relatively free graded algebras]{A model for the relatively free graded algebra of block-triangular matrices with entries from a graded PI-algebra}
\author{Lucio Centrone}\address{IMECC, Universidade Estadual de
Campinas, Rua S\'ergio Buarque de Holanda, 651, Cidade Universit\'aria ``Zeferino Vaz'', Distr. Bar\~ao Geraldo, Campinas, S\~ao Paulo, Brazil, CEP
13083-859}\email{centrone@ime.unicamp.br}
\author{Thiago Castilho de Mello}\thanks{T. C. de Mello was partially supported by grants from Capes (AUX-PE-PRODOC-2548/2010) and  from FAPESP  (No.
2012/16838-0)}\address{Instituto de Matematica e Estat\'{i}stica, Universidade de S\~{a}o Paulo Caixa Postal 66281, S\~{a}o Paulo, SP, 05315-970
Brasil}\email{tcmello@ime.unicamp.br}\keywords{}\subjclass[2000]{}
\begin{abstract} Let $G$ be a group and $A$ be a $G$-graded algebra satisfying a
polynomial identity. We construct a model for the relatively free $G$-graded algebra of the $G$-graded algebra of block-triangular matrices with entries
from $A$. We obtain, as an application, the ``factoring'' property for the $T_{\Z_2}$-ideals of block-triangular matrices with entries from the infinite
dimensional Grassmann algebra $E$ for some particular $\Z_2$-grading. \end{abstract}\maketitle

\section{introduction}
Let $F$ be a field and we denote by $\FF$ the free associative algebra freely generated by the set $X$ over $F$. We refer to elements of $\FF$ as polynomials in the
non-commutative variables of $X$. If $A$ is an $F$-algebra we say that $A$ satisfies a polynomial identity (or $A$ is a PI-algebra) if there exists an element
$f=f(x_1,\dots,x_n)\in \FF$ such that $f(a_1,\dots,a_n)=0$, for any $a_1,\dots,a_n\in A$. Such $f$ is called a polynomial identity of $A$. The set of all
polynomial identities of $A$, denoted by $T(A)$, is called the T-ideal of $A$. It is an ideal of $\FF$ which is invariant under endomorphisms of the algebra $\FF$.
The algebra $\FF/(\FF\cap T(A))$ is called the \textit{relatively free algebra} of $A$ (see the paper of Belov \cite{bel1} for a detailed account on relatively free algebras and their representability in any characteristic).

We recall that if $R$ is a block-triangular matrix algebra such that \[R=\left(
\begin{array}{cc}
A & M\\
0 & B
\end{array}
\right),\]where $A$ and $B$ are PI algebras and $M$ is an $A,B$-bimodule, then $R$ is a PI-algebra too. We say that $R$ (or also $T(R)$) satisfies the \textit{factoring property} if $T(R)=T(A)T(B)$, where the product of T-ideals is the usual product of ideals. We refer to the papers \cite{dre2} and \cite{dre1} by Drensky for a conjecture regarding the factoring property of T-ideals of block-triangular matrices with entries from minimal algebras already solved positively by Giambruno and Zaicev in \cite{giz3}.

In this paper we focus on the factoring property of block-triangular matrices in the graded case. Let $G$ be a group and let us consider a block
triangular matrix algebra \[R=\left(
\begin{array}{cc}
A & M\\
0 & B
\end{array}
\right),\]where $A$ and $B$ are PI $G$-graded algebras and $M$ is an $A,B$-bimodule. In \cite{dil1} Di Vincenzo and La Scala gave a proof for the Theorem of Lewin
in the graded case. Moreover, if $A$ is a $G$-graded subalgebra of $M_n(F)$ they introduce the notion of $G$-regularity for $A$. We recall that a map
$|\cdot|:\{1,\ldots,n\}\rightarrow G$ induces a $G$-grading over $M_n(F)$ such that for each $i,j$ the homogeneous degree of the matrix unit $e_{ij}$ is
$|j||i|^{-1}$. If $A$ is a graded subalgebra of $M_n(F)$, graded by a finite abelian group $G$, we have that $A$ is $G$\textit{-regular} if $|\cdot|$ is
surjective and its fibers are equipotent (see Theorem 5.4 of \cite{dil1}). They also proved that if $A\subseteq M_n(F)$ and $B\subseteq M_m(F)$ are matrix algebras graded
by a finite abelian group $G$, then $R$ has the factoring property for its ideal of graded polynomial identities if and only if one between $A$ and $B$ is
$G$-regular.

In the proof of the Theorem of Lewin, both in the ordinary and in the graded case, one prominent role is played by the relatively free (graded) algebra.
In \cite{pro1} Procesi showed the $k$-generated relatively free algebra of $M_n(F)$ is isomorphic to the $k$-generated algebra of generic matrices over the
polynomial ring, and such construction can be generalized for relatively free algebras of finite dimensional PI-algebras (see \cite{row1}). Moreover in \cite{ber1} Berele constructed the $k$-generated relatively free algebras of the minimal algebras $M_n(E)$ and $M_{a,b}(E)$.
It turned out that they are isomorphic to some $k$-generated subalgebra of generic matrices over the supercommutative polynomial algebra.

Let $A$ be an algebra over a field $F$. If $d_1,d_2,\dots,d_m$ are positive integers, we denote by $UT(d_1,\dots,d_m;A)$ the subalgebra of the matrix algebra $M_{d_1+\cdots+d_m}(A)$ consisting of matrices of the type

\[\left(
\begin{array}{ccccc}
A_{11} & A_{12}  & \ldots & A_{1m}\\
0 & A_{22}  & \ldots & A_{2m}\\
\vdots &  & \ddots & \vdots \\
0 & \ldots & \ldots  & A_{mm}
\end{array}
\right),\]
where $A_{ij}\in M_{d_i\times d_j}(A)$ for each $i,j$. One such algebra is called the algebra of \textit{block-triangular matrices} of size $d_1$, \dots,
$d_m$ over $A$. We may observe that if in addition $A$ is a PI-algebra, then $UT(d_1,\ldots,d_m;A)$ is a PI-algebra, too. In this paper we consider a PI $G$-graded algebra $A$ and we construct a model for the relatively free $G$-graded algebra of
$UT(d_1,\ldots,d_m;A)$. It turns out that the relatively free $G$-graded algebra of $UT(d_1,\ldots,d_m;A)$ is isomorphic to the algebra of generic matrices
with entries from the graded relatively free algebra of $A$. We use the model to prove that $UT(d_1,\ldots,d_m;A)$ has the factoring property when the relatively free graded
algebra of $A$ has a partially multiplicative basis. We also observe that the factoring property fails if we consider $UT(d_1,\ldots,d_m;E)$ with the grading
induced by the map $|\cdot|_{k^*}$ of the work of Di Vincenzo and Da Silva (see \cite{dit1}). The paper is organized as follows. Sections 2 is dedicated to the
main definitions concerning the graded polynomial identities. In Sections 3 we present the main tool of the paper: the graded version of the Theorem of Lewin
obtained by Di Vincenzo and La Scala in \cite{dil1}. In Section 4 we present the results regarding the $\Z_2$-graded identities of the infinite dimensional Grassmann
by Di Vincenzo and Da Silva (see \cite{dit1}). The results of Sections 4 will be used for applications of the main Theorem.
Section 5 is devoted to the construction of the model for the relatively free graded algebra of the upper triangular block matrices with entries from a $G$-graded
PI algebra. We obtain the factoring property for $UT(d_1,\ldots,d_m;A)$, where $A$ is a $G$-regular algebra having the same $T_G$ ideal of $M_n(F)$
when $G$ is a finite abelian group. In Section 6 we present the main theorem and, as a consequence, we obtain the factoring property for
$UT(d_1,\ldots,d_m;E)$ when $\Z_2$-graded by the $\Z_2$-grading induced by the natural grading over $E$. We also mention that the factoring property for
$UT(d_1,\ldots,d_m;E)$ has been proved by Berele and Regev in the ungraded case (see \cite{bere1}).

\section{Graded structures}
All fields we refer to are assumed to be of characteristic zero and all algebras we consider are associative and unitary.

Let $(G,\cdot)=\{g_1,\ldots,g_r\}$ be any finite group, and let $F$ be a field. If $A$ is an associative $F$-algebra, we say that $A$ is a
$G$\textit{-graded} algebra (or $G$-algebra) if there are subspaces $A^g$ for each $g\in G$ such that \[A=\bigoplus_{g\in G}A^g \ \textrm{and} \
A^gA^h\subseteq A^{gh}.\] If $0\neq a\in A^g$ we say that $a$ is \emph{homogeneous of $G$-degree g} or simply that $a$ has degree $g$, 
and we write $\deg(a)=g$. We consider the following subset of $G$ \[\support(A)=\{g\in G|A^g\neq0\}.\] The latter is called the \textit{support} of the $G$-graded algebra $A$.

We define a free object; let $\{X^{g}\mid g \in G\}$ be a family of disjoint countable sets. Put $X=\bigcup_{g\in G}X^{g}$ and denote by $\F$ the free
associative algebra freely generated by the set $X$. An indeterminate $x\in X$ is said to be of \emph{homogeneous $G$-degree} $g$, written $\deg(x)=g$,
if $x\in X^{g}$. We always write $x^{g}$ if $x\in X^{g}$. The homogeneous $G$-degree of a monomial $m=x_{i_1}x_{i_2}\cdots x_{i_k}$ is defined to be
$\deg(m)=\deg(x_{i_1})\cdot\deg(x_{i_2})\cdot\cdots\cdot\deg(x_{i_k})$. For every $g \in G$, we denote by $\F^g$ the subspace of $\F$ spanned by all the
monomials having homogeneous $G$-degree $g$. Notice that $\F^g\F^{g'}\subseteq \F^{gg'}$ for all $g,g' \in G$. Thus \[\F=\bigoplus_{g\in G}\F^g\] proves
$\F$ to be a $G$-graded algebra. The elements of the $G$-graded algebra $\F$ are referred to as $G$-graded polynomials or, simply, graded polynomials.
An ideal $I$ of $\F$ is said to be a $T_{G}$-ideal if it is invariant under all $F$-endomorphisms $\varphi:\F\rightarrow\F$ such that
$\varphi\left(\F^g\right)\subseteq\F^g$ for all $g\in G$. If $A$ is a $G$-graded algebra, a $G$-graded polynomial $\fx$ is said to be a \emph{graded
polynomial identity} of $A$ if $f(a_1,a_2,\cdots,a_t)=0$ for all $a_1,a_2,\cdots,a_t\in\bigcup_{g\in G}A^g$ such that $a_k\in A^{\deg(x_k)}$,
$k=1,\cdots,t$. If $A$ has a non-zero graded polynomial identity, we say that $A$ is a \emph{G-graded polynomial identity algebra} (GPI-algebra). We denote by
$T_G(A)$ the ideal of all graded polynomial identities of $A$. It is a $T_G$-ideal of $\F$. If $A$ is ungraded, i.e.,
graded by the trivial group, we talk about polynomial identities and T-ideal of $A$. We recall that if the group $G$ is finite and $A$ is a $G$-graded
PI-algebra, then it satisfies a polynomial identity (see \cite{bgz1}, \cite{bec1}). Moreover, we recall that if two $GPI$-algebras $A$ and $B$ satisfy the
same graded identities, i.e., $T_G(A)=T_G(B)$, then they satisfy the same identities, i.e., $T(A)=T(B)$.

When one deals with $T_G$-ideals, given a subset $Y\subseteq X$ one can talk about the least $T_G$-ideal of $\F$ which contains the set $Y$. Such
$T_G$-ideal will be denoted by $\langle Y\rangle ^{T_G}$ and will be called the $T_G$-ideal generated by $Y$. We say that elements of $\langle Y\rangle
^{T_G}$ are consequences of elements of $Y$, or simply that they follow from $Y$. Given a $G$-algebra $A$ one of the main problems in PI-theory is to find
a finite set $Y$ such that $T_g(A)=\langle Y\rangle ^{T_G}$. Such a $Y$ is called a basis for the $G$-graded polynomial identities of $A$.

We denote by $\U_G(A)$ the factor algebra $\dfrac{\F}{\F\cap T_G(A)}$ and we shall call it the \textit{relatively free $G$-algebra of $A$}. We observe that if $G=\{1_G\}$, we obtain the definition of the relatively free algebra.


\section{$\Z_2$-graded identities for the Grassmann algebra}
In this section we recall the main tools and definitions that are necessary for the study of graded polynomial identities of the Grassmann algebra. We
shall indicate the infinite dimensional Grassmann algebra by $E$.

The algebra $E$ can be constructed as follows. Let $\FF$ be the free algebra of countable rank on $X=\{x_1,x_2,\ldots\}$. If $I$ is the two-sided ideal of
$\FF$ generated by the set of polynomials $\{x_ix_j+x_jx_i|i,j\geq1\},$ then $E=\FF/I.$ If we write $e_i=x_i+I$ for $i=1,2,\ldots,$ then $E$ has the
following presentation:$$E=\langle1,e_1,e_2,\ldots|e_ie_j=-e_je_i,\text{\rm for all $i,j\geq1$}\rangle.$$ We say that the vector space $V$ generated by $X$ over $F$ is the generating vector space of $E$.
Moreover, the set
\[B=\{1,e_{i_1}\cdots e_{i_k}|1\leq i_1<\cdots<i_k\}\] is a basis of $E$ over $F$. Sometimes it is convenient to write $E$ in the form $E=E^{0}\oplus E^{1}$, where
\[E^{0}:=\text{\rm span}\{1,e_{i_1}\cdots e_{i_{2k}}|1\leq i_1<\cdots<i_{2k},k\geq0\},\]
\[E^{1}:=\text{\rm span}\{e_{i_1}\cdots e_{i_{2k+1}}|1\leq i_1<\cdots<i_{2k+1},k\geq0\}.\]
It is easily verified that the decomposition $E=E^{0}\oplus E^{1}$ is a $\Z_2$-grading of $E$ called the natural grading.  Notice that $E^{0}$ coincides
with the center of $E$. We give a look at the whole class of homogeneous $\Z_2$-grading of $E$. For more details we refer to the work of Di Vincenzo and Da
Silva (\cite{dit1}).

For a homogeneous $\Z_2$-grading of $E$ we mean any $\Z_2$-grading such that the generating vector space $V$ is a homogeneous subspace. This is equivalent
to consider a map $$\deg:V\rightarrow\Z_2.$$ If $w=e_{i_1}e_{i_2}\cdots e_{in}\in E$, then the set $\text{\rm Supp}(w):=\{e_{i_1},e_{i_2},\ldots,e_{i_n}\}$
is the support of $w$ and we define the $\Z_2$-grading of $w$ by $$\deg(e_{i_1}e_{i_2}\cdots e_{i_n})=\deg(e_{i_1})+\cdots+\deg(e_{i_n}).$$ If, for all
$e_i\in B$, one has $\deg(e_i)=1\in\Z_2,$ then we obtain the natural $\Z_2$-grading on $E$.

In this case, let $E^0$ be the homogeneous component of $\Z_2$-degree 0 and let $E^1$ be the component of degree 1. As we said above, $E^0=Z(E)$ is the center of $E$ and $ab + ba = 0$ for all $a, b\in E^1$. This means that $E$ satisfies the following graded polynomial identities:
$[y_1, y_2]$, $[y_1, z_1]$, $z_1z_2 + z_2z_1$.
Now, let us consider the $\Z_2$-gradings on $E$ induced by the maps  $\deg(\cdot)_{k*}$, $\deg(\cdot)_\infty$, and $\deg(\cdot)_k$, defined respectively by:

\[\deg(e_i)_{k*}=\left\{\begin{array}{ll}
\text{\rm 1 for $i = 1,\ldots, k$}\\
\text{\rm 0 otherwise},\end{array}\right.\]

\[\deg(e_i)_{\infty}=\left\{\begin{array}{ll}
\text{\rm 1 for $i$ odd}\\
\text{\rm 0 otherwise},\end{array}\right.\]

\[\deg(e_i)_k=\left\{\begin{array}{ll}
\text{\rm 0 for $i = 1,\ldots, k$}\\
\text{\rm 1 otherwise}.\end{array}\right.\]

By \cite{dit1}, we have the following result.

\begin{theorem}\label{grass}
Let $Y$ be a countable set of indeterminates of degree 0 and $Z$ be a countable set of indeterminates of degree 1 and put $X=Y\cup Z$. Then:
\begin{enumerate}
\item The $T_{\Z_2}$-ideal of $E$ graded by $\deg(\cdot)_\infty$ is generated by the polynomial \[[x_1,x_2,x_3].\]
\item The $T_{\Z_2}$-ideal of $E$ graded by $\deg(\cdot)_{k^*}$ is generated by the polynomials \[[x_1,x_2,x_3], z_1z_2\cdots z_{k+1}.\]
\item The $T_{\Z_2}$-ideal of $E$ graded by $\deg(\cdot)_{k}$ is generated by the polynomials \[[x_1,x_2,x_3],
 \text{\rm $[y_1,y_2]\cdots[y_{k-1},y_k][y_{k+1},x]$ (if $k$ is even)}, \]
\[\text{\rm $[y_1,y_2]\cdots[y_{k},y_{k+1}]$ (if $k$ is odd)},\]
\[ \text{\rm $g_{k-l+2}(z_1,\ldots,z_{k-l+2})[y_1,y_2]\cdots[y_{l-1},y_l]$ (if $l\leq k$)},\]
\[\text{\rm $[g_{k-l+2}(z_1,\ldots,z_{k-l+2}),y_1][y_2,y_3]\cdots[y_{l-1},y_l]$ (if $l\leq k,$ $l$ is odd)},\]
\[\text{\rm $g_{k-l+2}(z_1,\ldots,z_{k-l+2})[z,y_1][y_2,y_3]\cdots[y_{l-1},y_l]$ (if $l\leq k,$ $l$ is odd)},\]
\end{enumerate}
where the $g_m=g_m(z_1,\dots,z_m)$ are some polynomials in the odd variables $z_1,\dots,z_m$.
\end{theorem}

\begin{remark}\label{monomialidentity}
We note that among the three $\mathbb{Z}_2$-gradings defined above, the Grassmann algebra satisfies a graded monomial identity only if its $\Z_2$-grading is induced by
$\deg(\cdot)_{k^*}$.
\end{remark}

\begin{remark}\label{basis}
It is easily verified that a basis for the relatively free $\Z_2$-graded algebra of $E$ with the grading induced by $\deg(\cdot)_\infty$ is the following:
\[y_{i_1} y_{i_2}\cdots y_{i_n}z_{j_1}\cdots z_{j_m} [x_{l_1},x_{l_2}]\cdots[x_{l_{2s-1}},x_{l_{2s}}],\] where $n\geq 0$, $i_1\leq \cdots\leq  i_{n}$,
$m\geq 0$, $j_1\leq \cdots \leq j_m$, $s\geq 0$ and $l_1<\cdots <l_{2s}$.

On the other hand, a basis for the relatively free $\Z_2$-graded algebra of $E$ with the grading induced by $\deg(\cdot)_{k^*}$ is the following:
\[y_{i_1} y_{i_2}\cdots y_{i_n}z_{j_1}\cdots z_{j_m} [x_{l_1},x_{l_2}]\cdots[x_{l_{2s-1}},x_{l_{2s}}],\] where $n\geq 0$, $i_1\leq \cdots\leq  i_{n}$,
$0\leq m\leq k$, $j_1\leq \cdots \leq j_m$, $s\geq 0$ and $l_1<\cdots <l_{2s}$.
\end{remark}

\section{The graded Thorem of Lewin}
We resume the work of Di Vincenzo and La Scala (see \cite{dil1}) for the generalization of the Theorem of Lewin (see \cite{lew1}) in the graded case. In what follows the grading group is supposed to be finite.

Let $A$, $B$ be $G$-graded algebras and $M$ be an $A$,$B$-bimodule, then it is possible to consider the $G$-algebra \[R=\left(
\begin{array}{cc}
A & M\\
0 & B
\end{array}
\right).\]Let us call $\{x_i\}$'s the generators of the relatively free graded algebra of $A$ and $B$.
Then the following result holds.
\begin{proposition}\label{Lewin}Let $A$ and $B$ be GPI-algebras.
If $M$ contains a countable free set $\{u_i\}$ of homogeneous elements such that $\deg(x_i)=\deg(u_i)$ for any $i\geq1$, then
$T_G(R)=T_G(A)T_G(B)$.
\end{proposition}

We fix a map $|\cdot|:\{1,2,\cdots,n\}\rightarrow G$. Then $|\cdot|$ induces a grading on $M_n(F)$ by setting $|e_{ij}|= |j||i|^{-1}$, for all matrix units
$e_{ij}$. Indeed this is an elementary grading defined by $(|1|,\ldots, |n|)$. 
If we assume $A$ to be a matrix algebra, i.e., $A\subseteq M_n(F)$, with an elementary $G$-grading, the authors introduce the notion of $G$-regularity.
In the case $G$ is abelian, we have the following equivalence result.

\begin{theorem}
The $G$-algebra $A$ is $G$-regular if and only if the map $|\cdot|$ is surjective and all its fibers are equipotent.
\end{theorem}

Note that the $G$-regularity of $A$ is verified when the order of $G$
is exactly $n$ and the map $|\cdot|$ is bijective. This is the case, for instance, when we consider
the Di Vincenzo and Vasilovsky $\Z_n$-grading of $M_n(F)$ (see \cite{div1} and \cite{vas1}). Moreover, for the ordinary case, that is for $G=\{1_G\}$, the algebra $A$ is regular. We close the section with the main result of \cite{dil1}.

\begin{theorem}\label{gradedlewin}
Let $R$ be the $G$-graded block-triangular matrix algebra defined as above, where $A\subseteq M_n(F)$, $B\subseteq M_m(F)$ are $G$-algebras and $M$
is an $A$,$B$-bimodule. If one between $A$ and $B$ is $G$-regular, then the $T_G$-ideal $T_G(R)$ factorizes as: $T_G(R) = T_G(A)T_G(B)$.
\end{theorem}

\begin{corollary}\label{factoringmatrix}
Let \[R=\left(
\begin{array}{ccccc}
A_{11} & A_{12}  & \ldots & A_{1m}\\
0 & A_{22}  & \ldots & A_{2m}\\
\vdots &  \vdots& \ddots & \vdots \\
0 & 0 & \ldots  & A_{mm}
\end{array}
\right)\] be a $G$-graded subalgebra of some matrix algebra. Suppose that each $A_{ii}$ is a $G$-regular $G$-graded subalgebra of some $G$-graded matrix
algebra $M_{d_i}(F)$ and that $A_{ij}=M_{d_i\times d_j}(F)$, for each $i$ and $j$. Then
\[T_G(R)=T_G(A_{11})T_G(A_{22})\cdots T_G(A_{mm}).\]
\end{corollary}

\proof
We prove the statement by induction on $n$. If $n=2$, we are in the hypothesis of Theorem \ref{gradedlewin} and we are done. Suppose the assertion true for
$m-1$, where $m\geq 3$. We consider \[A=\left(
\begin{array}{ccccc}
A_{11} & A_{12}  & \ldots & A_{1m-1}\\
0 & A_{22} & \ldots  & A_{2m-1}\\
\vdots & & \ddots & \vdots \\
0 & \ldots & \ldots  & A_{m-1m-1}
\end{array}
\right),\] then \[R=\left(
\begin{array}{cc}
A & M\\
0 & A_{mm}
\end{array}\right).\] Due to the fact that $R$ is a $G$-algebra, we have $A$ is a $G$-algebra too. Now we are again in the hypothesis of Theorem
\ref{gradedlewin} and the proof follows.
\endproof

\section{A model for the relatively free graded algebra of block-triangular matrices}
In this section we construct a model for the relatively free graded algebra of $UT(d_1,\ldots,d_m;A)$, where $A$ is a GPI-algebra. We recall that in the finite dimensional case there is a standard way to construct such a model. For more details we refer to the book of Rowen \cite{row1}.

Let $A$ be a PI $G$-algebra over a field $F$. If $d_1,d_2,\dots,d_m$ are positive integers, we denote by $UT(d_1,\dots,d_m;A)$ the subalgebra of the matrix algebra $M_{d_1+\cdots+d_m}(A)$ consisting of matrices of the type

\[\left(
\begin{array}{ccccc}
A_{11} & A_{12}  & \ldots & A_{1m}\\
0 & A_{22}  & \ldots & A_{2m}\\
\vdots & \vdots & \ddots & \vdots \\
0 & 0 & \ldots  & A_{mm}
\end{array}
\right),\] where $A_{ij}\in M_{d_i\times d_j}(A)$ for each $i,j$. One such algebra is called the algebra of \textit{block-triangular matrices} of size
$d_1$, \dots, $d_m$ over $A$. In what follows we are going to construct a model for the relatively free graded algebra of $UT(d_1,\dots,d_m;A)$, where $A$
is any GPI-algebra.

We shall use the following notation: if $f(x_1,\ldots,x_n)$ is a graded polynomial of $\F$, we shall indicate by $\overline{x}$ the string of the
homogeneous indeterminates appearing in $f$, i.e., $\overline{x}=(x_1,\ldots,x_n)$, and we shall write $f(\overline{x})$ instead of $f(x_1,\ldots,x_n)$.
Moreover, if we are dealing with any graded substitution of the type $x_1^{\deg(x_1)}\mapsto a_1^{\deg{(x_1)}}$, we shall indicate the valuation of $f$ by
$f(\overline{a})$.

The model is based on the following construction. For each $k\in \mathbb{N}$ and $g\in G$, we define the matrix $\xi_k^{(g)}\in UT(d_1,\dots,d_m;\U_G(A))$
by
\[\xi^{(g)}_k=\left(
               \begin{array}{cccc}
                 B^{(g)}_{d_1\times d_1,k} & B^{(g)}_{d_1\times d_2,k} & \cdots & B^{(g)}_{d_1\times d_m,k} \\
                 0 & B^{(g)}_{d_2\times d_2,k} & \cdots & B^{(g)}_{d_2\times d_m,k} \\
                 \vdots & \vdots & \ddots & \vdots \\
                 0 & 0 & \cdots & B^{(g)}_{d_m\times d_m,k} \\
               \end{array}
             \right)\]

Here, each $B^{(g)}_{d_r\times d_s,k}$ is a $d_r\times d_s$ matrix whose entry $(i,j)$ is $x_{ij,k}^{(g)}+T_G(A) \in \U_G(A)$, with $d_1+\cdots+
d_{r-1}+1\leq i \leq d_1+\cdots+d_r$ and $d_1+\cdots+ d_{s-1}+1\leq j \leq d_1+\cdots+d_s$.

We denote by  $\U(d_1,\dots,d_m;A)$ the subalgebra of $UT(d_1,\dots,d_m;\U_G(A))$, generated by the matrices $\xi^{(g)}_k,\,k\in\mathbb{N}, g\in G$,
defined above (we omit an index $G$ in order to simplify the notation).

\begin{lemma}\label{Generic}
The algebra $\U(d_1,\dots,d_m;A)$ is a generic model for the relatively free graded algebra of $UT(d_1,\dots,d_m;A)$, i.e.,
\[\U(d_1,\dots,d_m;A)\cong \frac{\F}{T_G(UT(d_1,\dots,d_m;A))}.\]
\end{lemma}

\proof  Let $X=\bigcup_{g\in G}X^g$, where the union is disjoint and each $X^g$ is a countable set of homogeneous indeterminates. Define the homomorphism
\[\begin{array}{cccc} \varphi: & \F & \longrightarrow
& \U(d_1,\dots,d_m;A)\\ & x_k^{(g)} & \mapsto & \xi^{(g)}_k.\\
\end{array}\]
Of course $\varphi$ is a graded homomorphism onto $\U(d_1,\dots,d_m;A)$. We shall show \[\ker \varphi = T_G(UT(d_1,\dots,d_m;A)).\] First, observe that
$\ker \varphi \subseteq T_G(UT(d_1,\dots,d_m;A))$. Indeed, if $f=f(\overline{x})\in \ker \varphi$, then $\varphi(f)=f(\overline{\xi})=0$. Since $\U_G(A)$
is the relatively free algebra of $A$, we have that for any $\{a_{ij,k}^{(g)}\}\subseteq A$, there exists a graded homomorphism $\U_G(A)\longrightarrow A$,
such that $x_{ij}^{(g)}+T_G(A)\mapsto a^{(g)}_{ij,k}$. Since $f(\overline{\xi})=0$ in $\U(d_1,\dots,d_m;A)$, the image of any entry of the matrix
$f(\overline{\xi})$ under any homomorphisms $\U_G(A)\longrightarrow A$ is zero. Then $f\in T_G(UT(d_1,\dots,d_m;A))$.

In order to show the reverse inclusion, we take $f(\overline{x})\in \F$ a graded polynomial identity for $UT(d_1,\dots,d_m;A)$ and we consider the matrix
$M=f(\overline{\xi})$. Each entry $(r,s)$ of $M$ has the form $m_{rs}(x_{ij,k}^{(g)})+T_G(A)$, for some polynomials $m_{rs}\in \F$, with $1\leq r,s\leq
d_1+\cdots+d_m$. We claim that $M=0$. Indeed, since $f(\overline{u})=0$, for any $\overline{u}\in UT(d_1,\dots,d_m;A)$, we have that
$m_{rs}(a_{ij,k}^{(g)})=0$, for any $a_{ij,k}^{(g)}\in A$ and this shows that $m_{rs}\in T_G(A)$, for any $r$ and $s$, and hence $M=0$, which concludes the
lemma. \endproof

\begin{remark}
Notice that the above construction is still valid if we consider algebras over arbitrary fields.
\end{remark}

Now we present some applications of the above construction. In particular, we are interested in studying the ``factoring'' property of T-ideals.

\begin{lemma}\label{frac}
Let $A$ and $B$ be GPI-algebras such that $T_G(A)=T_G(B)$. Then \[T_G(UT(d_1,\dots,d_m;A))=T_G(UT(d_1,\dots,d_m;B)).\] As a consequence, for $n\in
\mathbb{N}$, we have $T_G(M_n(A))=T_G(M_n(B))$.
\end{lemma}

\proof Let $\xi^{(g)}_k$ be the generators of $\U(d_1,\dots, d_m;A)$ as in the previous theorem, and $\eta_k$ be the generators of $\U(d_1,\dots,d_m;B)$. Let
\[f(\overline{x})\in T_G(UT(d_1,\dots,d_m;A)),\] then we claim that $f\in T_G(UT(d_1,\dots,d_m;B))$. Let $m_{rs}(x_{ij,k}^{(g)}+T_G(A))$ be the entry $(r,s)$ of the matrix
$f(\overline{\xi})$, where $m_{rs}$ are polynomials in the variables $x_{ij,k}^{(g)}$. Since $f\in T_G(UT(d_1,\dots,d_n))$, by Lemma \ref{Generic} we have
$f(\overline{\xi})=0$, then $m_{rs}(x_{ij,k}^{(g)}+T_G(A))=0$, i.e., $m_{rs}\in T_G(A)=T_G(B)$. Hence, $m_{rs}(x_{ij,k}^{(g)}+T_G(B))=0$, for every $(r,s)$. Then we have
$f(\overline{\eta})=(m_{rs}(x_{ij,k}^{(g)}+T_G(B)))=0$, and $f\in T_G(UT(d_1,\dots,d_m;B))$. The other inclusion is analogous. \endproof


\begin{theorem}\label{main3}
Let $G$ be a finite abelian group and $A$ be a GPI-algebra such that $T_G(A)=T_G(M_n(F))$, where $M_n(F)$ is $G$-regular. Then for any set of positive
integers $\{d_1,\dots, d_m\}$, we have
\[T_G(UT(d_1,\dots,d_m;A))=T_G(M_{d_1}(A))T_G(M_{d_2}(A))\cdots T_G(M_{d_m}(A)).\]
\end{theorem}

\proof If $T_G(A)=T_G(M_n(F))$, by Lemma \ref{Generic} we have \[T_G(UT(d_1,\dots,d_m;A))=T_G(UT(d_1,\dots,d_m; M_n(F))).\] In light of the fact that
$M_n(F)$ is $G$-regular, the grading is induced by the $n$-tuple $\overline{g}=(g_1,\ldots,g_n)$. The $n$-tuple $\overline{g}$ gives rise to the map
$|\cdot|:\{1,\ldots,n\}\rightarrow G$ such that its fibers are equipotent. It is easy to see that $M_k(M_n(F))$ is $G$-graded isomorphic to $M_{nk}(F)$,
for any $k$, where the $G$-grading is induced by the $kn$-tuple \[\overline{g'}=(g_1,\ldots,g_n,g_1,\ldots,g_n,\ldots,g_1,\ldots,g_n).\] Hence
\[T_G(UT(d_1,\dots,d_m; M_n(F)))=T_G(UT(nd_1,\dots,nd_m;F))).\] Since each of the $M_{nd_{i}}(F)$'s is $G$-regular, it follows by Corollary
\ref{factoringmatrix} that the right-hand side of the above equation is equal to \[T_G(M_{nd_1}(F))\cdots T(M_{nd_m}(F)).\] By Lemma \ref{frac}, for each
$i$, we have
\[T_G(M_{nd_i}(F))=T_G(M_{d_i}(M_n(F)))=T_G(M_{d_i}(A)),\] and the result follows.

\endproof

%

\section{The factoring property}
In this section we shall use the model we constructed in the previous section in order to obtain a ``factoring'' theorem for block-triangular matrices with
entries from a $G$-graded algebra $A$, whose relatively free graded algebra has a partially multiplicative basis. As a consequence, we show that the
Grassmann algebra $E$, when $\Z_2$-graded by the grading inherited by the natural $\Z_2$-grading of $E$, has the factoring property.

Let $A$ be a GPI-algebra and consider the automorphism $\varphi$ of $\U_G(A)$ defined by $\varphi(x_{ij,k}^{(g)})=x_{i+1 j+1,k}^{(g)}$, where the
sum on the indexes is taken modulo $n$. Roughly speaking, we consider the automorphism of $\U_G(A)$, which sends each entry $(i,j)$ of a generic matrix to the
next variable in its diagonal.

The next lemma shows how the entries of the elements of $\U(n;A)$ behave with respect to this automorphism.
\begin{lemma}\label{independence}
Let $A$ be a GPI-algebra and $M=(m_{ij})\in \U(n;A)$. Then, for each $i$ and $j$, we have $\varphi(m_{ij})=m_{i+1 j+1}$.
\end{lemma}

\proof We use the same notations of Lemma \ref{Generic}. It is enough to prove the lemma for monomials and we do it by induction on the degree of such
monomials. We consider the monomial $f=f(B_1,\dots,B_r)$, where the $B_t$'s are the generic matrices in $\U(n;A)$. If $f$ has degree one, then $f$ is a
scalar multiple of one of such generic matrices and the result holds. Suppose the result is true for monomials of degree strictly less than $m$. If $f$ has
degree $m$, we write $f=g B_k$, where $g$ has degree $m-1$ and $g=g(B_1,\dots,B_r)=(p_{ij})$. If $B_k=(x_{ij,k}^{(g)})$, we have
$f=\left(\sum_{t=1}^np_{it}x_{tj,k}^{(g)}\right)$. The $(i,j)$ entry of $f$ is $ \sum_{t=1}^np_{it}x_{tj,k}^{(g)}$. Hence
$\varphi(\sum_{t=1}^np_{it}x_{tj,k}^{(g)})=\sum_{t=1}^n\varphi(p_{it})\varphi(x_{tj}^{(k)})=\sum_{t=1}^np_{i+1 t+1}x_{t+1 j+1,k}^{(g)}$ that is the $(i+1,
j+1)$ entry of $f$, and we are done. \endproof

As a consequence, the next result shows that in order to verify whether the elements of $\U_G(n;A)$ are linearly independent,  it is enough to verify the
linear independence of their $k$-th columns.

\begin{corollary}\label{linearindependence}
Let $A$ be a GPI-algebra, $\U_G(n;A)$ be the graded generic algebra of $M_n(A)$ and $\{f_1,\dots f_n\} \subseteq \U_G(n;A)$. If $f^k_i$ is the
$k$-th column of $f_i$ then the set of column vectors with entries from $\U_G(A)$, $\{f^k_1,\dots,f^k_n\}$ is linearly independent over $F$ if and only
if $\{f_1,\dots f_n\}$ is linearly independent over $F$.
\end{corollary}

\begin{remark}
With the above notations, one can see that for each $k$, the map
\[\pi_k: \sum_{i,j=1}^n a_{ij}e_{ij} \mapsto \sum_{i=1}^n a_{ik}e_{ik}\] is an injective map from $\U_G(n;A)$ into $M_n(\U_G(A))$.
\end{remark}

Before stating our main theorem, we introduce the concept of partially multiplicative basis.

\begin{definition}
If $A$ is a GPI-algebra, we say that the basis $B$ of $\U_G(A)$ is \emph{partially multiplicative}, if for any $S_1$, $S_2\subseteq B$ such that the elements of $S_1$ and the elements of $S_2$ are polynomials in disjoint sets of graded variables, the set $S_1S_2=\{s_1s_2\,|\, s_1\in S_1, s_2\in S_2\}$ is
linearly independent over $F$.
\end{definition}

\begin{theorem}\label{last}
Let $A$ be a GPI-algebra, graded by a group $G$ such that $\support(A)=G$. Suppose that its relatively free graded algebra, $\U_{G}(A)$, has a partially multiplicative basis. If
$d_1,\dots, d_m$ are positive integers, then
\[T_{G}(UT(d_1,\dots,d_m;A))=T_{G}(M_{d_1}(A))\cdots T_{G}(M_{d_m}(A)).\]
\end{theorem}

\proof We prove the theorem by induction on $m$. If $m=1$ the result is obvious. Suppose the assertion true for $m-1$, where $m>1$. We use the
same notation of the proof of Lemma \ref{Generic}. By Lemma \ref{Generic} we have
\[T_{G}(UT(d_1,\dots,d_m;A))=T_{G}(\U(d_1,\dots,d_m;A)).\] The
algebra $\U(d_1,\dots,d_m;A)$ is generated by the matrices
\[\xi_k^{(g)}=\left(
               \begin{array}{cccc}
                 B^{(g)}_{d_1\times d_1,k} & B^{(g)}_{d_1\times d_2,k} & \cdots & B^{(g)}_{d_1\times d_m,k} \\
                 0 & B^{(g)}_{d_2\times d_2,k} & \cdots & B^{(g)}_{d_2\times d_m,k} \\
                 \vdots & \vdots & \ddots & \vdots \\
                 0 & 0 & \cdots & B^{(g)}_{d_m\times d_m,k} \\
               \end{array}
             \right), \quad k\geq 1,\]where $g$ ranges over $G$. We consider now the algebra $C$ generated by the matrices \[\omega_k^{(g)}=\left(
               \begin{array}{cccc}
                 B^{(g)}_{d_1\times d_1,k} & B^{(g)}_{d_1\times d_2,k} & \cdots & B^{(g)}_{d_1\times d_{m-1},k} \\
                 0 & B^{(g)}_{d_2\times d_2,k} & \cdots & B^{(g)}_{d_2\times d_{m-1},k} \\
                 \vdots & \vdots & \ddots & \vdots \\
                 0 & 0 & \cdots & B^{(g)}_{d_{m-1}\times d_{m-1},k} \\
               \end{array}
             \right), \quad k\geq 1,\] and the algebra $B$ generated by the matrices $\eta_k^{(g)}=B_{d_m\times d_m,k}^{(g)}$, for $k\geq 1$ and $g\in G$.

We recall that by Lemma \ref{Generic}, the algebra $C$ is a relatively free $G$-graded algebra of $UT(d_1,\dots,d_{m-1};A)$ freely generated by the
matrices $\omega_k^{(g)}$ and $B$ is a relatively free algebra of $M_{d_m}(A)$, freely generated by the matrices $\eta_k^{(g)}$.

Let $d=d_1+\dots+d_{m-1}$, and let $M$ be the $(C,B)$-bimodule generated by the $d\times d_m$ matrices \[\mu_k^{(g)}=\left(
                                                                   \begin{array}{c}
                                                                     B^{(g)}_{d_1\times d_m,k} \\
                                                                     B^{(g)}_{d_2\times d_m,k} \\
                                                                     \vdots \\
                                                                     B^{(g)}_{d_{m-1}\times d_m,k} \\
                                                                   \end{array}
                                                                 \right), \quad k\geq 1,\] with action given by the usual product of matrices.

By Proposition \ref{Lewin} (the graded theorem of Lewin), if $M$ is freely generated by the $\mu_k^{(g)}$, we have
\[T_{G}(\U(d_1,\dots,d_m;A))=T_{G}(\U(d_1,\dots,d_{m-1};A))T_{G}(\U(d_m;A))\] and the proof follows by induction.

Due to the fact that $\support(A)=G$, we have that the $\mu_k^{(g)}$'s are non-zero. In order to show that the bimodule $M$ is freely generated by the $\mu_k^{(g)}$, we take $\{\mu_1,\dots,\mu_m\}$ a finite subset of $\{\mu_k^{(g)}\}$, and
we suppose $\{f_1,\dots,f_n\}\subseteq B$ to be a linearly independent set over $F$ and $\{h_{ij}\}\subseteq C$. We need to show that if
$\sum_{i=1}^n\sum_{j=1}^m h_{ij}\mu_jf_i=0$, then each $h_{ij}$ is zero in $C$. We have $\sum_{j=1}^m\left(\sum_{i=1}^n h_{ij}\mu_jf_i\right)=0$, and
each $\mu_j$ depends on disjoint sets of variables. Then for each $j\in \{1,\dots, m\}$, $\sum_{i=1}^n h_{ij}\mu_jf_i$ depends on a different set of
variables. Hence, if $\sum_{i=1}^n(\sum_{j=1}^m h_{ij}\mu_jf_i)=0$, we also have that $\sum_{i=1}^n h_{ij}\mu_jf_i=0$, for each $j$, then it is enough
to prove the assertion for only one of the $\mu_j$, say $\mu$, i.e., we need to prove that if $\{f_1,\dots,f_n\}\subseteq B$ is linearly independent over $F$
and $\{h_1,\dots, h_n\}\subseteq C$ is such that $\sum_{i=1}^nh_i\mu f_i=0$, then each $h_i$ is zero.

We put  $d=d_1+\dots+d_{m-1}$ and rename the variables $x_{ij,k}^{(g)}$,  $1\leq i\leq d$, $d+1\leq j\leq d+d_m$ from $\mu_k^{(g)}$ by $y_{ij,k}^{(g)}$,
$1\leq i\leq d$, $1\leq j\leq d_m$ and the variables $x_{ij,k}^{(g)}$, $d+1 \leq i,j\leq d+d_m$ from $\eta_k^{(g)}$ by $z_{ij,k}^{(g)}$, $1\leq i,j\leq
d_m$.

For each $t$, $h_t=h_t(\omega_1^{(g_1)},\dots,\omega_{k_t}^{(g_t)})$ is a matrix of the form $(h_{rs,t})$, where each
$h_{rs,t}=h_{rs,t}(x_{ij,k}^{(g)})$ is a polynomial in the variables $x_{ij,k}^{(g)}$, for $1\leq i,j\leq d$ and $k\geq 1$. The matrix $\mu$ has the
form $(y_{pq})$, $1\leq p\leq d$ and $1\leq q\leq d_m$.

With this notation, multiplying matrices we have
\[h_t \mu=\left(\sum_{l=1}^d h_{rl,t}y_{ls}\right)_{1\leq r\leq d,\; 1\leq s\leq d_m}.\]

If $f_t=f_t(\eta_1^{(g_1)},\dots,\eta_{n_t}^{(g_t)})$, we can write it as a matrix of the form $(f_{pq,t})$, where $f_{pq,t}$ is a polynomial in the
variables $z_{ij,k}^{(g)}$, $1\leq i,j\leq d_m$ and $k\geq 1$. Then for each $t$, we have
\[h_t \mu f_t=\left(\sum_{k=1}^{d_m}\sum_{l=1}^d h_{rl,t}y_{lk}f_{ks,t}\right)_{1\leq r\leq d,\; 1\leq s\leq d_m}.\]

If $\sum_{t=1}^n h_t \mu f_t=0$, we obtain
\[0=\sum_{t=1}^n h_t \mu f_t=\left(\sum_{t=1}^n\sum_{k=1}^{d_m}\sum_{l=1}^d h_{rl,t}y_{lk}f_{ks,t}\right)_{1\leq r\leq d,\; 1\leq s\leq d_m}.\]

Hence each entry of the above matrix is zero, i.e., for $1\leq r\leq d$ and $1\leq s\leq d_m$, \[\sum_{t=1}^n\sum_{k=1}^{d_m}\sum_{l=1}^d
h_{rl,t}y_{lk}f_{ks,t}=0.\]

Since the polynomials $h_{lr,t}$ depend only on the variables $x$'s, and the polynomials $f_{ks,t}$ depend only on the variables $z$'s, comparing the
degree of the variables $y_{kl}$ in the above sum, we have that for each $k$ and $l$
\[\sum_{t=1}^n h_{rl,t}y_{lk}f_{ks,t}=0.\]

For $r=s=l=1$ we have the system of equations

\[\displaystyle\sum_{t=1}^n h_{11,t}y_{1s}f_{s1,t}=0 \quad 1\leq s\leq d.\]


Since $\{f_1,\dots,f_n\}$ is linearly independent over $F$, Corollary \ref{linearindependence} implies that the set
\[\left\{\left(\begin{array}{c}
                                                                                                  f_{11,t} \\
                                                                                                      \vdots \\
                                                                                                  f_{d1,t}  \\
                                                                                                \end{array}\right)
,\; 1\leq t\leq n\right\}\subseteq M_{d\times 1}(\U_{G}(A))\] is also linearly independent over $F$.

For each $t\in \{1,\dots, n\}$, and $s\in \{1,\dots, d\}$ write
\[f_{s1,t}=\sum_{i=1}^m\alpha_{s1,t}^{(i)}b_i \text{ and }\]
\[h_{11,t}=\sum_{j=1}^r\beta_{t}^{(j)}c_j,\]
where $b_i$ and $c_j$ are elements of a partially multiplicative basis of $\U_G(A)$. Since the $f's$ and the $h's$ depend on disjoint sets of variables, the
set $\{c_jy_{1s}b_i,\ |\ i,j\} $ is linearly independent in $\U_G(A)$.

This means that in the above system of equations, we have

\[\sum_{t=1}^n\sum_{j=1}^r\sum_{i=1}^m\alpha_{s1,t}^{(i)}\beta_{t}^{(j)}c_jy_{1s}b_i=0,\quad 1\leq s\leq d,\]

that is equivalent to

\[\sum_{j=1}^r\sum_{i=1}^m\left(\sum_{t=1}^n\alpha_{s1,t}^{(i)}\beta_{t}^{(j)}\right)c_jy_{1s}b_i=0, \quad 1\leq s\leq d.\]

Since the set $\{c_jy_{1s}b_i\}$ is linearly independent over $F$, the above system is equivalent to

\begin{equation}\label{system}
\sum_{t=1}^n\alpha_{s1,t}^{(i)}\beta_{t}^{(j)}=0\quad \text{ for each } 1\leq s\leq d, \; 1\leq i\leq m,\; 1\leq j\leq r.
\end{equation}

It is enough to prove that the only solution of the above system is $\beta_t^{(j)}=0$, for $1\leq t\leq n$ and $1\leq j\leq r$.

Fix $1\leq j\leq m$; since the set \[\left\{\left(\begin{array}{c}
                                                                                                  f_{11,t} \\
                                                                                                      \vdots \\
                                                                                                  f_{d1,t}  \\
                                                                                                \end{array}\right)
,\; 1\leq t\leq n\right\}\subseteq M_{d\times 1}(\U_{G}(A))\] is linearly independent over $F$,
\[\sum_{t=1}^n \beta_t^{(j)}\left(\begin{array}{c}
                                                                                                  f_{11,t} \\
                                                                                                      \vdots \\
                                                                                                  f_{d1,t}  \\
                                                                                                \end{array}\right)=0\] if and only if
                                                                                                $\beta_t^{(j)}=0$ for every $t$.
We show now that the above equation is equivalent to the system of equations (\ref{system}).

Substituting the expressions for the $f_{s1,t}$, the above equation is equivalent to
\[\left(\begin{array}{c} \beta_1^{(j)}(\alpha_{11,1}^{(1)}b_1+\cdots +\alpha_{11,1}^{(m)}b_m)+\cdots+\beta_n^{(j)}(\alpha_{11,n}^{(1)}b_1+\cdots+ \alpha_{11,n}^{(m)}b_m)\\
                         \beta_1^{(j)}(\alpha_{21,1}^{(1)}b_1+\cdots +\alpha_{21,1}^{(m)}b_m)+\cdots+\beta_n^{(j)}(\alpha_{21,1}^{(1)}b_1+\cdots +\alpha_{21,1}^{(m)}b_m)\\
                         \vdots \\
                         \beta_1^{(j)}(\alpha_{d1,1}^{(1)}b_1+\cdots +\alpha_{d1,1}^{(m)}b_m)+\cdots+\beta_n^{(j)}(\alpha_{d1,1}^{(1)}b_1+\cdots +\alpha_{d1,1}^{(m)}b_m) \end{array}\right)=0.\]
Factoring the $b_i$, we obtain
\[\left(\begin{array}{c} (\beta_1^{(j)}\alpha_{11,1}^{(1)}+\cdots+\beta_n^{(j)}\alpha_{11,n}^{(1)})b_1+\cdots+(\beta_1^{(j)}\alpha_{11,1}^{(m)}+\cdots+\beta_n^{(j)}\alpha_{11,n}^{(m)})b_m\\
                         (\beta_1^{(j)}\alpha_{21,1}^{(1)}+\cdots+\beta_n^{(j)}\alpha_{21,n}^{(1)})b_1+\cdots+(\beta_1^{(j)}\alpha_{21,1}^{(m)}+\cdots+\beta_n^{(j)}\alpha_{21,n}^{(m)})b_m\\
                         \vdots \\
                         (\beta_1^{(j)}\alpha_{d1,1}^{(1)}+\cdots+\beta_n^{(j)}\alpha_{d1,n}^{(1)})b_1+\cdots+(\beta_1^{(j)}\alpha_{d1,1}^{(m)}+\cdots+\beta_n^{(j)}\alpha_{d1,n}^{(m)})b_m \end{array}\right)=0.\]
Due to the fact that the $b_i$ are linearly independent over $F$, we conclude that
\[\sum_{t=1}^n\beta_t^{(j)}\alpha_{s1,t}^{(i)}=0, \text{ for each } 1\leq s\leq d,\; 1\leq i\leq m,\;1\leq j\leq r\]
that is exactly the system of equations (\ref{system}), i.e., we have shown that the only solution to the system (\ref{system}) is the trivial
one. As a consequence, we have that $h_{11,t}=0$, for every $1\leq t\leq n$.

Similar calculations show that $h_{pq,t}=0$, for every $1\leq p,q\leq d$ and $1\leq t\leq n$, which concludes the proof of the theorem.
\endproof%

\section{Conclusions}
We observe that the model we presented in this paper transfers the information regarding the graded identities of a graded algebra $A$ into the ideal of graded
polynomial identities of a block-triangular matrix $UT(d_1,\ldots,d_m;A)$. We may observe that the map $\deg(\cdot)_\infty$ induces a grading over $E$ such
that its relatively free $\Z_2$-graded algebra has a partially multiplicative basis. We can state the following.

\begin{theorem}\label{last2}
Let $E$ be the Grassmann algebra endowed with the $\Z_2$-grading induced by the map $\deg(\cdot)_\infty$. If $d_1,\dots, d_m$ are positive integers, then
\[T_{\Z_2}(UT(d_1,\dots,d_m;E))=T_{\Z_2}(M_{d_1}(E))\cdots T_{\Z_2}(M_{d_m}(E)).\]
\end{theorem}

As we already noted in Remark \ref{monomialidentity}, the only homogeneous $\Z_2$-grading over $E$ that gives rise to a monomial identity is the one
induced by the map $\deg(\cdot)_{k^*}$. Indeed, the relatively free $\Z_2$-graded algebra does not have a partially multiplicative basis. In fact, the Theorems
\ref{last} and \ref{last2} cannot be generalized. We have the following.

\begin{proposition}
Let $R$ be the $\Z_2$-graded algebra \[R:=\left(
\begin{array}{cc}
E & E\\
0 & E
\end{array}
\right),\]where the $\Z_2$-grading is induced by that of $E$. If $E$ is graded by $\deg(\cdot)_{k^*}$,  then \[T_{\Z_2}(E)T_{\Z_2}(E)\varsubsetneq
T_{\Z_2}(R).\]
\end{proposition}
\proof It is easy to see that $z_1\cdots z_{k+1}$ is a graded identity of $R$ but it is clearly not a consequence  of the product $T_{\Z_2}(E)T_{\Z_2}(E)$
because $T_{\Z_2}(E)$ contains the identity $z_1\cdots z_{k+1}$ (see Theorem \ref{grass}).
\endproof

It is well known that the (ungraded) Grassmann algebra, $E$, has a partially multiplicative basis for its relatively free-graded algebra. As a consequence,
Theorem \ref{last} gives that $UT(d_1,\dots,d_m;E)$ has the factoring property for its T-ideal.

A natural question arises: ``what are the PI-algebras such that their relatively free algebras do not have a partially multiplicative basis?"

\end{document}